\documentclass[12pt]{amsart}
\theoremstyle{definition}

\usepackage{url,upref,amscd,amsmath,amsthm}
\usepackage{times}
\usepackage{fancyhdr}
\usepackage[psamsfonts]{amssymb}
\usepackage{cite}
\begin{document}
\title[MOR cryptosystem & extra-special $p$-groups]{The MOR cryptosystem and extra-special $p$-groups}
\author{Ayan Mahalanobis}\address{Indian Institute of Science
  Education and Research Pune, Pashan Pune-411021, India}
\today
\email{ayan.mahalanobis@gmail.com}
\keywords{MOR cryptosystem, extra-special $p$-groups, the discrete
  logarithm problem}
\thanks{This research is supported by a NBHM research grant}
\begin{abstract} This paper studies the MOR cryptosystem, using the
  automorphism group of the extra-special $p$-group of exponent $p$,
  for an odd prime $p$. Similar results can be obtained for
  extra-special $p$-groups of exponent $p^2$ and for the even prime.
\end{abstract}
\maketitle
\section{Introduction}
In this paper, we study the MOR cryptosystem with extra-special $p$
groups. Similar studies were done, using the group of unitriangular
matrices~\cite{ayan3} and the group of unimodular
matrices~\cite{ayan4}. The group of unitriangular
matrices and the group of unimodular matrices are both matrix
groups. There are many ways to represent a group -- natural
representations, like a matrix representation or permutation
representation, or a more abstract representation in the form of
generators and relations, commonly known as a \emph{finite presentation}.
In this paper, we shift our study of the MOR cryptosystem,
from the matrix representation of a group to a finite
  presentation. We show that using finite presentation, in the form of
generators and relations, one can build a \textbf{secure} MOR
cryptosystem.

In a MOR cryptosystem, one works with the \emph{discrete logarithm problem}
in the automorphism group. On one hand, this is not a major change, because the
discrete logarithm problem works in a group and the automorphisms form
a group. On the other hand, an automorphism group arises from any
algebraic structure, like a graph, vector space, etc. So the MOR cryptosystem
can be seen, as the one, that liberates the discrete logarithm problem
from groups to other algebraic structures.

\section{The MOR cryptosystem}
In this section we describe the MOR cryptosystem~\cite{crypto2001,asiacrypt2004} as
automorphisms of a finite group $G$, however it can be generalized to
other finitely generated algebraic structures easily. A description and a critical analysis of the MOR cryptosystem is in
\cite{ayan3} and the references there. 
\subsection{Description of the MOR cryptosystem}\label{MOR}
Let $G=\langle g_1,g_2,\ldots,g_\tau\rangle$, $\tau\in\mathbb{N}$
be a finite group and $\phi$ a non-trivial automorphism of 
$G$. Alice's keys are as follows: 
\begin{description}\label{keyex}
\item[Private Key] $m$, $m\in\mathbb{N}$.
\item[Public Key] $\left\{\phi(g_i)\right\}_{i=1}^\tau$ and $\left\{\phi^m(g_i)\right\}_{i=1}^\tau$.
\end{description}
\paragraph{\textbf{Encryption}}
\begin{description}
\item[a] To send a message (plaintext) $a\in G$ Bob computes $\phi^r$
  and $\phi^{mr}$ for a random $r\in\mathbb{N}$.
\item[b] The ciphertext is $\left(\left\{\phi^r(g_i)\right\}_{i=1}^\tau,\phi^{mr}(a)\right)$.
\end{description}
\paragraph{\textbf{Decryption}}
\begin{description}
\item[a] Alice knows $m$, so if she receives the ciphertext
  $\left(\phi^r,\phi^{mr}(a)\right)$, she computes $\phi^{mr}$ from $\phi^r$ and
  then $\phi^{-mr}$ and then computes $a$ from $\phi^{mr}(a)$.
\end{description}
Alice knows the order of the automorphism $\phi$, she can use
the identity $\phi^{t-1}=\phi^{-1}$ whenever $\phi^t=1$ to compute
$\phi^{-mr}$.
\section{Notations and definitions}
All definitions are standard and so are the notations.
\begin{description}
\item The exponent of a finite group $G$ is the least common multiple of all
possible orders of elements in $G$. For a finite $p$-group, it is
the largest order of an element in $G$.
\item The center of a group $G$, denoted by Z$(G)$, is the set of all
  elements in $G$ that commute with every element of $G$. It is known
  that Z$(G)$ is \emph{characteristic}.
\item For a group $G$, $G^\prime$ is the commutator of $G$ and
  $\Phi(G)$ is the Frattini subgroup of $G$, see~\cite[Page 2]{green}
  for details.
\end{description}
\section{The description and some properties of extra-special $p$-groups}
For a given prime $p$, all groups of order $p^2$ are abelian. So the
first non-abelian group $G$ is of order $p^3$. There is a complete
classification of groups of order $p^3$. For $p=2$, there are two groups
of of order $8$, the dihedral group
$D_8$, and the quaternion group $Q_8$.
\subsection{Groups of order $p^3$, for an odd prime $p$}
For a odd prime $p$, there are two non-isomorphic
classes~\cite[Section 4.13]{suzuki2} of
non-abelian groups of order $p^3$:
\begin{eqnarray}
M:=\langle x,y\;|\;x^p=1=y^p; [x,y]=z\in\text{Z}(M); z^p=1\rangle\\
N:=\langle x,y\;|\;y^p=1; [x,y]=x^p=z\in\text{Z}(N); z^p=1\rangle
\end{eqnarray}
Both of these groups are 2-generator $p$-groups, the first one has exponent
$p$ and the second one has exponent $p^2$. In this paper we study the
MOR cryptosystem using $M$, similar study can be done with $N$ and with
the $D_8$ and $Q_8$, with similar conclusions. Let $\phi$ be an
automorphism of $M$, then $\phi$ 
can be written as
\begin{eqnarray}
\phi(x)=x^{m_1}x^{n_1}z^{l_1}\\
\phi(y)=x^{m_2}x^{n_2}z^{l_2}.
\end{eqnarray}
Then $[\phi(x),\phi(y)]=z^{\det(T)}$, where 
$T=\begin{pmatrix}
m_1&n_1\\
m_2&n_2
\end{pmatrix}$. This shows that $\det(T)\neq 0\mod p$. Notice that
$\dfrac{M}{\Phi(M)}\cong\mathbb{Z}_p\times\mathbb{Z}_p$, and $M$ is
\emph{extra-special}, hence the group of inner automorphisms of $M$, denoted
by $I$, is isomorphic to $\mathbb{Z}_p\times\mathbb{Z}_p$. This gives
the following exact sequence:
\[
\begin{CD}
0 @>>>\mathbb{Z}_p\times\mathbb{Z}_p @>>>\text{Aut}(M)
@>>>\text{GL}(2,p) @>>> 1
\end{CD}
\]
There are two kinds of automorphisms of $M$, one that is trivial on
Z$(M)$ and the other that is not. Since any automorphism of the center
of $M$ can be extended to an automorphism of $M$, the automorphism
that acts non-trivially on the center are generated by
\begin{equation}\label{one} 
x\mapsto x, \;\;\; y\mapsto y^\theta
\end{equation}
where $\theta$ is primitive
mod $p$. If we denote the automorphisms that are trivial on the center
by $H$, then there is an exact sequence of the form
\[
\begin{CD}
0 @>>>\mathbb{Z}_p\times\mathbb{Z}_p @>>>H
@>>>\text{SL}(2,p) @>>> 1
\end{CD}
\]
Since for $M$, the central and the inner automorphisms are identical,
the inner automorphisms are of the form $x\mapsto xz^{d_1}$, $y\mapsto yz^{d_2}$,
where $0\leq d_1, d_2<p$.

Hence we have shown that any automorphism $\phi$ of $M$ is a
composition of automorphisms, (\ref{one}), inner automorphism and
an element from SL$(2,p)$.

It is not hard to see that if $\phi$ is given by
\begin{eqnarray*}
\phi(x)=x^{m_1}y^{n_1}z^{l_1}\\
\phi(y)=x^{m_2}y^{n_2}z^{l_2}
\end{eqnarray*}
and $\phi^m$ is given by
\begin{eqnarray*}
\phi^m(x)=x^{{m_1}^\prime}y^{{n_1}^\prime}z^{{l_1}^\prime}\\
\phi^m(y)=x^{{m_2}^\prime}y^{{n_2}^\prime}z^{{l_2}^\prime}
\end{eqnarray*}
then
\[\begin{pmatrix} m_1 & n_1\\ m_2 & n_2 \end{pmatrix}^m
= \begin{pmatrix}m_1^\prime & n_1^\prime\\ m_2^\prime &
  n_2^\prime\end{pmatrix}.\]

So the discrete logarithm problem in the automorphism $\langle\phi\rangle$ is converted to
the discrete logarithm problem in GL$(2,p)$. Once can use $m_i$ and
$n_i$, $i=1,2$ in $\phi$,
such that, the matrix $T$ is in SL$(2,p)$. The best
algorithm to solve the discrete logarithm problem in matrices is the Menezes-Wu algorithm~\cite{menezes1}. That algorithm finds the
eigenvalues of the matrix and the eigenvalues of the power of that matrix, and then try to
solve the discrete logarithm problem in those eigenvalues. So if
the characteristic polynomial corresponding to the matrix of $\phi$ is
irreducible then the complexity to solve the discrete logarithm
problem in $\phi$ and $\phi^m$ is identical to solving the discrete
logarithm problem in $\mathbb{F}_{p^2}$. 
\subsection{Extra-special p-groups of exponent $p$} An extra-special
group $P$ is a $p$-group, in which the center Z$(P)$, the commutator
$P^\prime$, and the Frattini subgroup $\Phi(P)$ are equal and
cyclic of order $p$~\cite[Definition 4.14]{suzuki2}. The two most important extra-special p-groups are $M$ and $N$
above. Extra-special $p$-groups are well studied and their
automorphism groups was described by Winter~\cite{winter}. We don't
want to redo all the work done by Winter but refer an interested
reader to his paper~\cite{winter}.

Let $P$ be the \emph{iterative central
  product}~\cite[Section2.2]{green} of $M$ with itself $r$ times. As
we know  $M$ is a group of order $p^3$ and exponent $p$. This makes
$P$ an extra-special $p$-group of exponent $p$.
The finite presentation for the group $P$ is the following~\cite[Page 33]{green}:
\[
P=\langle x_1,\ldots,x_r,y_1,\ldots,y_r\;|\; [x_i,y_j]=1,
i\neq j;\; [x_i,y_i]=z\in\text{Z}(P)\rangle
\]
each of $x_i,y_i$ and $z$ is of order $p$. 

One can define a non-degenerate, bilinear alternating form, $B$, on
$\dfrac{P}{\Phi(P)}$ as a vector space over $\mathbb{Z}_p$~\cite[Page 33]{green}. Let $x,y\in P$, and $\overline{x},\overline{y}$ be their
image in $\dfrac{P}{\Phi(P)}$. Then $B\left(\overline{x},\overline{y}\right)=c$,
where $[x,y]=z^c$.

Description of the automorphisms of $P$ involves three steps.
\begin{description}
\item[A] Find all automorphisms that are non-trivial on the center.
\item[B] Prove that an automorphism preserves the bilinear form if and
  only if it acts trivially on the center. Let $H$ be the subgroup of
  the automorphism group that acts trivially on the center.
\item[C] Prove that $H/I\cong\text{Sp}(2r,p)$. Where $I$ is the subgroup
  of inner automorphisms of $P$ and $\text{Sp}(2r,p)$ is the symplicitic group
  on the vector space $\dfrac{P}{\Phi(P)}$ over $\mathbb{Z}_p$,
  defined by the bilinear form $B$.
\end{description}
We briefly sketch the proof of the above three assertions, for
details, see~\cite{winter}. It is known that for an extra-special
$p$-group the inner automorphisms are identical to the central
automorphisms. Hence the inner automorphisms are given by 
\[x_i\mapsto x_iz^{d_i},\;\;y_i\mapsto y_iz^{d_i^\prime}\] where $0\leq
d_i,d_i^\prime<p$. Clearly there are $p^{2n}$ inner automorphisms of $P$.
\paragraph{(A)}
The automorphisms that doesn't act trivially on Z$(P)$ are given by
powers of $z\mapsto
z^\theta$, where $\theta$ is a primitive element mod $p$. Notice that
Z$(P)$ is a cyclic group of order $p$. Hence these automorphisms can be
defined by:
\begin{equation}
\theta:\;x_i\mapsto x_i, \;\; y_i\mapsto y_i^\theta
\end{equation}
where $\theta$ is primitive mod $p$. Clearly, $\theta$ is of order
$p-1$.
\paragraph{(B-C)} Corresponding to an automorphism $\phi$ of $P$, one can
trivially define an automorphism $\overline{\phi}$ on
$\dfrac{P}{\Phi(P)}$. Then the automorphism $\overline{\phi}$ preserves the
bilinear form $B$ if and only if $\phi$ acts trivially on $Z(P)$. This
follows from the equation
\[
\left[\phi(x),\phi(y)\right]=B\left(\overline{\phi(x)},\overline{\phi(y)}\right)=B(\overline{x},\overline{y})=[x,y].
\]
Hence there is a homomorphism $\tau:H\rightarrow\text{Sp}(2r,p)$. It is easy to
see that the kernel is the set of inner automorphisms $I$. This proves
that $H/I\cong\text{Sp}(2r,p)$.

By an argument identical to the MOR cryptosystem in $M$, one can
reduce the discrete logarithm problem in the extra-special $p$-group $P$ to
that of a discrete logarithm problem in Sp$(2r,p)$. The discrete
logarithm problem in Sp$(2r,p)$, in the best case scenario (irreducible
characteristic polynomial), embeds into an discrete logarithm problem
in $\mathbb{F}_{p^{2r}}$.
\section{Conclusion}
The discrete logarithm problem is the backbone of many modern day
public key cryptosystems and key exchanges. A MOR cryptosystem
generalizes the central idea of the discrete logarithm problem from
a group to any finitely generated algebraic structure.

It was an open question, if one can build a secure MOR cryptosystem using the
finite presentation of a group. We have shown that the answer is yes.

The situation with other extra-special $p$-groups is almost identical.
\begin{small}
\bibliography{paper}
\bibliographystyle{amsplain}
\end{small} 
\end{document}